\documentclass[a4paper,12pt]{article} 

\usepackage{fancybox} 
\usepackage{textcomp}
\usepackage{pifont}
\usepackage{siunitx}
\usepackage[stable]{footmisc}
\usepackage{graphicx}
\usepackage{caption}
\graphicspath{{pictures/}}
\DeclareGraphicsExtensions{.pdf,.png,.jpg}
\usepackage{comment}
\usepackage[russian]{babel}
\usepackage{fancybox} 
\usepackage{amssymb} 
\usepackage{amsmath} 
\usepackage[utf8]{inputenc}
\usepackage{amsfonts}
\usepackage{amsthm}
\usepackage{color}
\usepackage{mathrsfs} 
\usepackage{dsfont}
\usepackage{titling}
\usepackage{hyperref}
\hypersetup{
	colorlinks=true,
	linkcolor=blue,
	filecolor=magenta,      
	urlcolor=cyan,
}
\usepackage{mathtools}                                                      
\mathtoolsset{showonlyrefs=true}      
\usepackage{float}
\usepackage{color}

\newtheorem{theorem}{Теорема}[section]
\newtheorem{proposition}[theorem]{Предложение}
\newtheorem{corollary}[theorem]{Следствие}
\newtheorem{lemma}[theorem]{Лемма}

\newtheorem*{remark}{Замечание}

\def \R {{\mathbb {R}}}
\def \N {{\mathbb {N}}}

\def \H {{\mathbb {H}}}

\def \phi {{\varphi}}

\def \tilde {\widetilde}

\thanksmarkseries{arabic}
\numberwithin{equation}{section}

\title{\textbf{Асимптотический вариант метода параметрикс для цепей Маркова, сходящихся к диффузиям}}
\author{\textbf{И.Биттер}\thanks{Международная лаборатория стохастического анализа и его приложений, НИУ ВШЭ,
		Покровский бульвар, 11, г.Москва, Российская Федерация. ilya.bitter@yandex.ru}, \textbf{В. Конаков}\thanks{Международная лаборатория стохастического анализа и его приложений, НИУ ВШЭ,
		Покровский бульвар, 11, г.Москва, Российская Федерация. kv24@mail.ru}}

\begin{document}
	\maketitle

	\begin{abstract}
	В работе приводится обобщение локальной предельной теоремы о сходимости неоднородных цепей Маркова к диффузионному пределу на случай, когда соответствующие коэффициенты процессов удовлетворяют слабым условиям регулярности и совпадают лишь асимптотически. В частности, рассматриваемые нами коэффициенты сноса могут быть неограниченными с не более чем линейным ростом, а оценки отражают перенос терминального состояния неограниченным трендом через соответствующий детерминированный поток. Наш подход основан на изучении равномерного расстояния между переходными плотностями заданной неоднородной цепи Маркова и предельного диффузионного процесса, а оценка скорости сходимости получена с использованием классической локальной предельной теоремы и оценок устойчивости типа параметрикса.
	\end{abstract}
	
	\tableofcontents	
	\section{Введение} \label{Intro}
	\subsection{Постановка задачи}
{\allowdisplaybreaks	
	Для $n > n_0 \geq 1$ рассмотрим неоднородную цепь Маркова $X^n_{t_k}$, заданную на решётке $t_k = k/n$ с $k = 0 \cdots n$ и принимающую значения в $\mathbb{R}^d$. Динамика цепи $X^n$ имеет следующий вид:
	\begin{equation}\label{МС_0}
		X^n_{t_{i+1}} = X^n_{t_i} + \frac{1}{n}b^n(t_i, x) + \frac{1}{\sqrt{n}}\xi^n_{i+1}, \quad i = 1 \cdots n - 1,
	\end{equation}
где функция $b^n: [0, 1] \times \mathbb{R}^ d \rightarrow \mathbb{R}^d$, а семейство ошибок $\xi^n$ удовлетворяет следующему марковскому условию:
\begin{equation}\label{distr_1}
	\mathcal{L}\left(\xi^n_{i+1} | X^n_{t_i} = x_i, \cdots\right) = q^n_{t_i, x_i}(\cdot).
\end{equation}
Вероятностные распределения, отвечающие плотностям $q^n_{t_i, x_i}(\cdot)$, имеют нулевое среднее при любых значениях $n, t_i, x_i$, а соответствующие условные ковариационные матрицы определяются соотношением
\begin{equation}\label{cov_1}
	\int_{\mathbb{R}^d} z_i z_j q^n_{t,x}(z) dz = a^n_{ij}(t,x).
\end{equation}
Функция $a^n$ отображает произведение $[0,1] \times \mathbb{R}^d$ в пространство вещественнозначных положительно определенных матриц $d \times d$, удовлетворяющих условию равномерной эллиптичности. В свою очередь, коэффициенты сноса $b^n$ также удовлетворяют некоторым слабым условиям регулярности, которые будут уточнены в дальнейшем, а также имеют не более чем линейный рост. Переходную плотность процесса \eqref{МС_0} будем обозначать $p_n$.

Теперь рассмотрим $d$-мерный невырожденный диффузионный процесс 
\begin{equation}\label{SDE_0}
dX_t = b(t, X_t)dt + \sigma(t, X_t)dW_t, \quad t \in [0,1],
\end{equation} 
где $W_t$ обозначает $d$ - мерное броуновское движение на некотором фильтрованном полном вероятностном пространстве
$(\Omega,\mathcal F,(\mathcal F_t)_{t \geq 0},\mathbb P) $  с обычными условиями, а коэффициенты предполагаются измеримыми и регулярными. Также, матрица $a(t, x):=\sigma \sigma^{*}(t, x)$ является равномерно эллиптичной. Когда оба коэффициента $b, \sigma$ ограничены и непрерывны по Гёльдеру, то хорошо известно, что уравнение \eqref{SDE_0} допускает единственное слабое решение \cite{fried, kalash, KKM}. Кроме того, комбинируя классический метод параметрикса \cite{deck, KKM} и метод цепочек \cite{bass}, можно показать \cite{Ar1, Ar2}, что переходная плотность $p$ процесса \eqref{SDE_0} удовлетворяет двусторонней оценке с константой $C \geq 1$
$$C^{-1} g_{\lambda^{-1}}(t, x-y) \leq p(0, x, t, y) \leq C g_{\lambda}(t, x-y),$$ где $$g_{\lambda}(t, x)=t^{-d / 2} \exp \left(-\frac{\lambda|x|^{2}}{t}\right), \lambda \in (0,1], \quad t>0.$$

Когда снос неограничен и нелинеен, известно меньше результатов. Чтобы получить верхнюю оценку, нам нужно контролировать члены
ряда параметрикса, а в случае неограниченного сноса это становится деликатной проблемой, и необходимо вводить прямой поток, соответствующий переносу начального состояния или, что эквивалентно, обратный поток, соответствующий переносу конечного состояния \cite{DM, KMM, MPZ}. Для неограниченного, не более чем линейно растущего сноса верхняя граница для переходной плотности может быть получена с использованием метода усечения, введенного в \cite{DM}, см. также \cite{MPZ}. 

При рассматриваемых условиях и равномерной сходимости $b^n \rightarrow b, a^n \rightarrow a$ при $n \to \infty$ из, например, Теоремы 1 в \cite{Sk1} следует, что имеет место сходимость по распределению цепи \eqref{МС_0} к \eqref{SDE_0}. Слабая сходимость распределений дискретных по времени марковских процессов к диффузиям широко изучена. Классическая литература в этой области включает результаты Скорохода \cite{Sk2}, Струка и Варадана \cite{SV79}. Эти результаты получены вероятностными методами. Однако, впоследствии весьма распространенным стало использование аналитического подхода к рассмотрению сходимости переходных плотностей марковских цепей к диффузиям, см. \cite{KM, KKM}. А именно, применение метода параметрикса для параболических УРЧП и модификация этого метода для дискретных по времени марковских цепей позволили количественно оценить слабую сходимость, упомянутую выше. Этот подход также может быть использован для доказательства локальных предельных теорем для алгоритмов стохастической аппроксимации, известных как процедуры Роббинса-Монро \cite{KM0}. Более подробную информацию о~применении метода параметрикса можно найти в \cite{DM, MPZ}.

Отметим также, что к оценке скорости сходимости моделей \eqref{МС_0} и \eqref{SDE_0} приводит концепция индуцированных порядковых статистик. Асимптотическая теория индуцированных порядковых статистик обсуждается в \cite{B25,David,David2,Dav}.

Основной задачей данной работы является расширение результата, полученного в \cite{KM} на случай, когда пары коэффициентов $(b, a)$ и $(b^n, a^n)$ совпадают лишь асимптотически, а $b, b^n$ неограничены и допускают не более чем линейный рост. 

Статья организована следующим образом. В следующем подразделе мы вводим наши предположения и формулируем
основной результат. Раздел 2 содержит введение в некоторые важные факты о методе параметрикса в форме МакКина-Зингера и
детерминированном переносящем потоке, связанном со сносом соответствующей диффузии или цепи Маркова. В разделе 3 мы выводим наши основные результаты. Дополнительные доказательства технических лемм приведены в Приложении. 

В дальнейшем будем обозначать $\langle \cdot,\cdot\rangle$ и $|\cdot| $ евклидово скалярное произведение и норму на $\R^d$. Также $D^\nu_x = \prod\limits_{i = 1}^d D^{\nu_i}_{x_i}$ будет означать дифференцирование относительно мультииндекса $\nu=(\nu_1,\cdots,\nu_d)\in \N^d $, для которого  $|\nu|=\sum_{i=1}^d \nu_i$. Обозначение $C$ соответствует положительной константе, которая зависит только от введённых предположений; в контексте приводимых доказательств её значение может изменяться от строки к строке. 

\subsection{Предположения и основной результат}\label{results}
\begin{itemize}
	\item[\textbf{(A1)}] (\textbf{Равномерная эллиптичность.})
	Ковариационные матрицы процессов \eqref{МС_0} и \eqref{SDE_0} удовлетворяют условию равномерной эллиптичности при всех $(t,x) \in [0,1] \times \mathbb{R}^d$ с $\Lambda \geq 1$:
	\begin{align*}
		&\Lambda^{-1}|w|^{2} \leq \langle a(t, x) w, w\rangle  \leq \Lambda|w|^{2}, \\
		&\Lambda^{-1}|w|^{2} \leq \langle a^n(t, x) w, w \rangle \leq \Lambda|w|^{2}, \quad  \forall w \in \mathbb{R}^{d}.
	\end{align*}
	
	\item[\textbf{(A2)}] (\textbf{Регулярность.})
	Коэффициенты сноса и ковариационные матрицы в уравнениях \eqref{МС_0} и \eqref{SDE_0} удовлетворяют слабым условиям  регулярности, а именно: для всех пар $(t,x), (s,y) \in [0,1] \times \mathbb{R}^d$ существуют константы $A, B$ (вообще говоря, не зависящие от $n$) такие, что $B < n_0$ и
	\begin{align*}
		&\left|a(t, x) - a(s, y)\right| + \left|a^n(t, x) - a^n(s, y)\right| \leq A\left(\left|x-y\right|^{\gamma} + |t - s|^{\alpha}\right), \\
		&\left|b(t, x) - b(s, y)\right| + \left|b^n(t, x) - b^n(s, y)\right| \leq B\left(\left|x-y\right| + |t - s|^{\beta}\right), \quad \gamma, \beta, \alpha \leq 1.
	\end{align*}
Заметим, что при таких предположениях на коэффициенты процессов \eqref{МС_0} и \eqref{SDE_0} обратные переносящие потоки $\theta^n_{t_i, t_j}(y)$ и $\theta_{t_i, t_j}(y)$ могут быть определены как единственные решения разностного и дифференциального уравнений, соответственно:
\begin{equation}
	\begin{cases}
		\theta^n_{t_{i+1}, t_j}(y) = \theta^n_{t_i, t_j}(y) + \frac{1}{n}b_n(t_i, \theta^n_{t_i, t_j}(y)), \\
		\theta^n_{t_j, t_j}(y) = y.
	\end{cases}
	\begin{cases}
		\frac{d}{du} \theta_{u, s}(y) = b(u, \theta_{u, s}(y)), \\
		\theta_{s, s}(y) = y.
	\end{cases}
\end{equation}  
\item[\textbf{(A3)}] (\textbf{Неограниченность сноса.})
Коэффициенты сноса ограничены в нуле, то есть существует $K > 0$ такое, что при всех $t \in [0,1]$
\begin{equation*}
	\left|b(t, \textbf{0})\right| + \left|b^n(t, \textbf{0})\right| \leq K.
\end{equation*}
Принимая во внимание липшицевость по пространственной переменной, несложно понять, что указанное условие эквивалентно не более чем линейному росту трендов:
\begin{equation*}
	|b(t, x)| + |b^n(t, x)| \leq K\left(1 + |x|\right) \quad \forall (t, x) \in [0,1] \times \mathbb{R}^d.
\end{equation*}
Кроме того, компоненты коэффициентов сноса и ковариационной матрицы дважды непрерывно дифференцируемы по пространственной переменной, а производные удовлетворяют условиям \textbf{(A2)} и \textbf{(A3)}.
	\item[\textbf{(A4)}] (\textbf{Регулярность семейства плотностей.})
	Семейство плотностей \eqref{distr_1}, соответствующих распределениям ошибок в \eqref{МС_0}, равномерно по $n$  полиномиально убывает на бесконечности вместе со своими производными, а также удовлетворяет условию регулярности с некоторой константой $C > 0$:
	\begin{align*}
	&\left|D_z^{\nu} q^n_{t,x}(z)\right| \leq C Q_S(z) \quad \forall (t,x) \in [0,1] \times \mathbb{R}^d, z \in \mathbb{R}^d,  \\
	&\left|q^n_{t,x}(z) - q^n_{t,y}(z)\right| \leq \frac{C}{\sqrt{n}}\left|x - y\right| Q_S(z) \quad \forall (x,y,z) \in (\mathbb{R}^d)^3
	\end{align*}
	с полиномиальным ядром $Q_S(z) = \frac{C_S}{(1 + |z|)^S}$ для некоторого $S > 2d + 6$ и $|\nu| \leq 4$, а константа $C_S$ выбрана для нормализации указанного выражения. Масштабированный аналог $\mathcal{Q}_S(\frac{z}{\sqrt{t}}) = t^{-d/2} Q_S(\frac{z}{\sqrt{t}})$.
	\item[\textbf{(A5)}] (\textbf{Асимптотическая устойчивость процедуры.})	
	Пара коэффициентов $(b^n, a^n)$ равномерно сходится к соответствующим коэффициентам уравнения \eqref{SDE_0}:
	\begin{align*}
		&\Delta_{b} = \sup\limits_{(t,x) \in [0,1] \times \mathbb{R}^d} \left|b(t,x) - b^n(t,x)\right| \underset{n \rightarrow \infty}{\longrightarrow} 0,\\
			&\Delta_{a} = \sup\limits_{(t,x) \in [0,1] \times \mathbb{R}^d} \left|a(t,x) - a^n(t,x)\right| \underset{n \rightarrow \infty}{\longrightarrow} 0.	
			\end{align*}
\end{itemize}	

Приведем теперь основной результат работы. Он состоит в оценке равномерного расстояния между переходными плотностями процессов \eqref{МС_0} и \eqref{SDE_0}, существование которых следует из используемых предположений. 

\begin{theorem}
	Для всякого $n$ и пар $(t_i, x), (t_j, y)$ существует константа $C > 0$, зависящая от параметров в предположениях такая, что справедлива верхняя оценка скорости сходимости переходной плотности $p_n$ марковской цепи, определяемой уравнением \eqref{МС_0} к переходной плотности $p$ процесса \eqref{SDE_0}:
	\begin{align}\label{Main_Thm}
		&\left|p(t_i, t_j, x, y) - p_n(t_i, t_j, x, y)\right| \leq \\ &\leq  C\ln(e(j-i))(1+|x|^{S-d-2} + |y|^{S-d-2}) \cdot \Delta^n \cdot \left(\mathcal{Q}_{S-d-6}\left(\frac{x - \theta^n_{t_i, t_j}(y)}{\sqrt{t_j - t_i}}\right) +  \mathcal{Q}_{S-d-6}\left(\frac{x - \theta_{t_i, t_j}(y)}{\sqrt{t_j - t_i}}\right)\right),
	\end{align}
где коэффициент $\Delta^n = \dfrac{1}{n^{\min(\gamma/2, \alpha, \beta)}} + \Delta_{b} + \Delta_{a}$.
\end{theorem}
\section{Разложение переходных плотностей диффузий и цепей Маркова в ряд параметрикса}
\subsection{Метод параметрикса в форме МакКина-Зингера}
Предположим, что выполняются условия \textbf{(A1)}, \textbf{(A2)} и \textbf{(A3)}. Предположения влекут существование по меньшей мере дважды непрерывно дифференцируемой переходной плотности процесса \eqref{SDE_0}, и, следовательно, позволят нам применить технику PDE, а именно, рассмотреть указанную переходную плотность как фундаментальное решение соответствующих уравнений Колмогорова \cite{DM, KKM}. Наряду с процессом \eqref{SDE_0} введем его "замороженный аналог" $\tilde{X}_u$ с динамикой 
\begin{equation}\label{SDE_1}
	d\tilde{X}_u = b(u, \theta_{u, s}(y))du + \sigma(u, \theta_{u, s}(y))dW_u,
\end{equation}
где $0 \leq u \leq s \leq 1$, $y \in \mathbb{R}^d$, и, как было указано ранее, обратный переносящий поток $\theta_{u, s}(y)$ является единственным и определенным глобально на $[0,1]$ решением задачи Коши $\frac{d}{du} \theta_{u, s}(y) = b(u, \theta_{u, s}(y))$ с терминальным условием $\theta_{s, s}(y) = y$. 

Переходная плотность процесса \eqref{SDE_1} является гауссовской с вектором средних и ковариационной матрицей, имеющими следующий вид при $t < s$:
\begin{align*}
	\tilde{\theta}_{t,s}(x) = x + \int\limits_{t}^{s} b(u, \theta_{u, s}(y))du = x + y - \theta_{t, s}(y), \quad
	\tilde{\mathcal{C}}_{t,s}=\int\limits_t^s\sigma\sigma^*(u,\theta_{u,s}(y))du.
\end{align*}
Инфинитезимальные операторы процессов \eqref{SDE_0} и \eqref{SDE_1}, соответственно, в момент времени $u \in [0,1]$ и для $\varphi \in C_{0}^{2}\left(\mathbb{R}^{d}, \mathbb{R}\right), z \in \mathbb{R}^{d}$ задаются соотношениями 
\begin{align*}
	L_{u} \varphi(z, y)=\frac{1}{2} \operatorname{Tr}\left(\sigma \sigma^{*}(u, z) D_{z}^{2} \varphi(z, y)\right)+\left\langle b(u, z), D_{z} \varphi(z, y)\right\rangle, \\
	\tilde{L}_{u} \varphi(z, y)=\frac{1}{2} \operatorname{Tr}\left(\sigma \sigma^{*}(u, \theta_{u, s}(y)) D_{z}^{2} \varphi(z, y)\right)+\left\langle b(u, \theta_{u, s}(y)), D_{z} \varphi(z, y)\right\rangle.
\end{align*}
Отметим, что коэффициенты генератора процесса \eqref{SDE_1} "заморожены" в некоторой (терминальной) точке $(s,y) \in [u,1] \times \mathbb{R}^d$.

Введем \textit{ядро параметрикса} 
\begin{equation*}
	H(t,s,x,y) = \left(L_t - \tilde{L}_t\right)\tilde{p}(t,s,x,y),
\end{equation*}
и определим операцию свертки $\otimes$ соотношением 
\begin{equation*}
	f \otimes g(t, s, x, y)=\int_t^s du \int_{\mathbb{R}^d} dz f(t, u, x, z) g(u, s, z, y).
\end{equation*}
Применяя процедуру сглаживания коэффициентов процесса \eqref{SDE_0}, мы приходим к разложению переходной плотности $p$ в бесконечный ряд
\begin{equation}\label{decomp_cont}
	p(t, s, x, y)=\tilde{p}(t, s, x, y)+\sum_{r=1}^{\infty} \tilde{p} \otimes H^r(t, s, x, y),
\end{equation}
где $H^1 = H$, а при $r > 1$ индуктивно $H^r = H^{r-1} \otimes H$. Подробное обоснование данного разложения см. в \cite{MPZ}. Основное свойство ряда \eqref{decomp_cont} состоит в том, что он позволяет приблизить, вообще говоря, негауссовскую плотность $p$ гауссовской $\tilde{p}$ и ее производными.

Приведем теперь аналогичное разложение для переходного ядра $p_n$ процесса \eqref{МС_0}. Для этого определим "замороженную" дискретную цепь Маркова $\tilde{X}^n_{t_{k}}$ с динамикой 
\begin{equation*}\label{MC_1}
	\tilde{X}^n_{t_{i+1}} = \tilde{X}^n_{t_i} + \frac{1}{n}b_n(t_i, \theta^n_{t_i, t_j}(y)) + \frac{1}{\sqrt{n}}\tilde{\xi}^n_{i+1},
\end{equation*}
где $\theta^n_{t_{i+1}, t_j}(y) = \theta^n_{t_i, t_j}(y) + \frac{1}{n}b_n(t_i, \theta^n_{t_i, t_j}(y))$ и $\theta^n_{t_j, t_j}(y) = y$, $i = 0 \cdots j - 1$, а закон распределения ошибок $\tilde{\xi}^n_{i}$ удовлетворяет марковскому соотношению
\begin{equation}\label{distr_2}
	\mathcal{L}\left(\tilde{\xi}^n_{i+1} | \tilde{X}^n_{t_i} = x_i, \cdots\right) = q^n_{t_i, \theta^n_{t_i, t_j}(y)}(\cdot).
\end{equation}

Инфинитезимальные операторы процессов \eqref{SDE_0} и \eqref{SDE_1} для $t_i = i/n$ и $\varphi \in C_{0}^{2}\left(\mathbb{R}^{d}, \mathbb{R}\right), x \in \mathbb{R}^{d}$ имеют вид 
\begin{align*}
	{L}_{t_i}^n \varphi(x):=n \mathbb{E}\left[\varphi\left({X}_{t_{i+1}}^{n}\right)|{X}^n_{t_i} = x\right] -n\varphi(x), \\
	\tilde{L}_{t_i}^n \varphi(x):=n \mathbb{E}\left[\varphi\left(\tilde{X}_{t_{i+1}}^{n}\right)|\tilde{X}^n_{t_i} = x\right] -n\varphi(x),
\end{align*}
а дискретная операция свертки $\otimes_n$ определена как 
\begin{equation*}
	f \otimes_n g\left(t_i, t_j, x, y\right)=\sum_{k=0}^{j-i-1} \frac{1}{n} \int_{\mathbb{R}^d} f\left(t_i, t_{i+k}, x, z\right) g\left(t_{i+k}, t_j, z, y\right) d z.
\end{equation*}
\textit{Дискретное ядро параметрикса} удовлетворяет 
\begin{equation*}
	H^n\left(t_i, t_j, x, y\right):=\left(L_{t_i}^n-\tilde{L}_{t_i}^{n}\right) \tilde{p}_n\left(t_{i}, t_j, x, y\right),
\end{equation*}
и, как следствие, справедливо разложение 
\begin{equation}\label{decomp_discr}
	p_n\left(t_i, t_j, x, y\right)=\sum_{r=0}^{j-i} \tilde{p}_n \otimes_n H^{n,r}\left(t_i, t_j, x, y\right),
\end{equation}
где $H^{n,1} = H^n$, а при $r > 1$ $H^{n,r} = H^{n, r-1} \otimes_n H^n$. Подробности приведены в \cite{KM0, KKM}.
\subsection{Сходимость переходной плотности диффузии к своему ряду параметрикса}

Заметим, что знак равенства в разложении \eqref{decomp_cont} носит формальный характер до того момента, пока не доказана абсолютная сходимость ряда параметрикса. 

\begin{lemma}
	Для любых $0 \leq t_i < t_j \leq 1$ и $x,y \in \mathbb{R}^d$, а также мультииндекса $\nu$ с $|\nu| \leq 4$ существует константа $C > 0$ такая, что справедлива оценка
	\begin{align}\label{froz_dens_bound}
		\left|D_x^{\nu}\tilde{p}\right|(t_i, t_j, x, y) &+  \left|D_x^{\nu}\tilde{p}_n\right|(t_i, t_j, x, y) \leq \\ &\leq \frac{C}{(t_j - t_i)^{|\nu|/2}}\left(\mathcal{Q}_{S-d-1-|\nu|}\left(\frac{x - \theta_{t_i, t_j}(y)}{\sqrt{t_j - t_i}}\right) + \mathcal{Q}_{S-d-1-|\nu|}\left(\frac{x - \theta^n_{t_i, t_j}(y)}{\sqrt{t_j - t_i}}\right)\right).
	\end{align}
\end{lemma}
\proof
Неравенство $$\left|D_x^{\nu}\tilde{p}\right|(t_i, t_j, x, y) \leq \frac{C}{(t_j - t_i)^{|\nu|/2}}\left(\mathcal{Q}_{S-d-1-|\nu|}\left(\frac{x - \theta_{t_i, t_j}(y)}{\sqrt{t_j - t_i}}\right)\right)$$ может быть получено прямыми вычислениями, учитывая, что $$|x|^{\varepsilon}e^{-|x|^2} \leq C e^{-c|x|^2}$$ для всякого $\varepsilon > 0$ при некотором $0 < c < 1$ и $C > 0$, а также существует $C > 0$ такое, что $e^{-|x|^2} \leq \frac{C}{1 +|x|^S}$. Для оценки второго слагаемого в левой части \eqref{froz_dens_bound} применим классическую локальную предельную теорему (см. \cite{BR}). Отметим, что $\tilde{p}_n(t_i, t_j, x, y)$ является  плотностью суммы \begin{equation*}
	x + \sum_{k=i}^{j-1}\frac{1}{n}b_n(t_k, \theta^n_{t_k, t_j}(y)) + \sum_{k = i}^{j-1}\frac{1}{\sqrt{n}}\tilde{\xi}^n_{k+1}
\end{equation*}
в точке $\theta^n_{t_i, t_j}(y) - x$.
Пусть $f^n_{i,j}$ - плотность нормированной суммы 
	$S_{i,j}^n = V^{-1/2}_{i,j}\sum_{k = i}^{j-1}\frac{1}{\sqrt{n}}\tilde{\xi}^n_{k+1},$
где матрица $V_{i,j} = \sum_{k = i}^{j-1}\frac{1}{n}a^n(t_{k}, \theta_{t_{k}, t_j}(y))$. Очевидно, эта матрица также удовлетворяет условию равномерной эллиптичности, и имеет место \begin{equation*}
	C^{-1}(t_j - t_i)^{-d/2} \leq \det  V^{-1/2}_{i,j} \leq C(t_j - t_i)^{-d/2}
\end{equation*}
с константой $C \geq 1$ ввиду ограниченности компонент матрицы $a^n$. Таким образом, \begin{equation*}
	\tilde{p}_n(t_i, t_j, x, y) = \det V^{-1/2}_{i,j} f^n_{i,j}(\theta^n_{t_i, t_j}(y) - x).
\end{equation*}
Аргументы, аналогичные доказательству Леммы 3.7 из \cite{KM}, показывают, что для плотности $f^n_{i,j}$ применима Теорема 19.3 из \cite{BR}, вследствие чего мы получаем утверждение \eqref{froz_dens_bound} с $|\nu| = 0$. Для доказательства утверждения в общем случае применим Теорему 19.3 из \cite{BR} для $D^{\nu}_x \tilde{p}_n$ с очевидными модификациями (см. \cite{KM, KKM}).
\endproof
\begin{remark}
	Верхняя оценка для гауссовской переходной плотности "замороженной" диффузии может быть записана в виде \begin{equation}\label{upp_bound}
		\left|D^{\nu}_x \tilde{p}\right|(t, s, x, y) \leq \frac{C}{(s-t)^{(|\nu|+d)/2}}\exp\left(-\dfrac{|\theta_{t, s}(y) - x|^2}{C(s-t)}\right), \quad t<s
	\end{equation}
(см. \cite{BK21, DM}). В свою очередь, правая часть указанного выражения при $|\nu| = 0$  является нижней оценкой для переходной плотности $\bar{p}$ некоторого вспомогательного диффузионного процесса к постоянным коэффициентом диффузии \cite{MPZ}:
\begin{equation*}\label{aux_diff}
	d\bar{X}_t = b\left(t, \bar{X}_t\right) d t+\lambda I d W_{t}
	\end{equation*}
с достаточно большим коэффициентом $\lambda$. Следовательно, \begin{equation}\label{maj_dens_bound}
\left|D^{\nu}_x \tilde{p}\right|(t, s, x, y) \leq \frac{C}{(s-t)^{|\nu|/2}} \bar{p}(t, s, x, y).
\end{equation}
Иногда в ходе доказательств мы будем автоматически обновлять  мажорирующую плотность $\bar{p}(t, s, x, y)$ , имея в виду следующую последовательность неравенств для параметра $\varepsilon > 0$:
\begin{align*}
	\dfrac{\left|\theta_{t,s}(y) - x\right|^{\varepsilon}}{(s - t)^{\varepsilon/2}} \bar{p}(t,s,x,y) \leq \dfrac{\left|\theta_{t,s}(y) - x\right|^{\varepsilon}}{(s - t)^{\varepsilon/2}}\frac{C_1}{(s-t)^{d/2}}\exp\left(-\frac{|\theta_{t,s}(y)-x|^2}{C_1(s-t)}\right) \leq \\ \leq \frac{C_2}{(s-t)^{d/2}}\exp\left(-\frac{|\theta_{t,s}(y)-x|^2}{C_2(s-t)}\right) \leq C\bar{p}(t,s,x,y).
\end{align*}
\end{remark}
\begin{lemma}
	Для любых $0 \leq t < s \leq 1$ и $x,y \in \mathbb{R}^d$,  существует константа $C > 0$ такая, что справедливо 
	\begin{equation}\label{ker_bound}
		|H|(t, s, x, y) \leq \dfrac{C}{(s-t)^{1-\gamma/2}}\bar{p}(t, s, x, y).
	\end{equation}
\end{lemma}
\proof Доказательство следует из оценки \eqref{maj_dens_bound}, предположения \textbf{(A2)} и замечания к Лемме 2.1.
\endproof
\begin{corollary}
	Правая часть разложения в ряд параметрикса сходится абсолютно и равномерно на $\left([0,1] \times \mathbb{R}\right)^2$.
\end{corollary}
\proof Используя определение операции свёртки $\otimes$, оценки \eqref{maj_dens_bound}, \eqref{ker_bound} и уравнение Колмогорова-Чепмэна для подходящей мажорирующей плотности $\bar{p}$, приходим к оценке 
\begin{equation}\label{conv_convergence}
	\left|\tilde{p} \otimes H^{n}\right|(t,s,x,y)	\leq C^{n + 1} \dfrac{\Gamma\left(\frac{\gamma}{2}\right)^{n}}{\Gamma(1 + \frac{n\gamma}{2})} (s-t)^{\frac{n\gamma}{2}} \cdot \bar{p}(t, s, x, y).
\end{equation}
Асимптотическое поведение гамма-функции влечет абсолютную сходимость ряда в правой части \eqref{decomp_cont}.
\endproof
\subsection{Некоторые свойства переносящих потоков}

В данном разделе мы перечислим некоторые свойства дискретного и непрерывного переносящих потоков $\theta^n_{t_i, t_j}(y)$ и $\theta_{t_i, t_j}(y)$ (см. предположение \textbf{(A2)}). Приводимые свойства будут важны для дальнейшего анализа.

\begin{proposition}
	\begin{itemize}
		Зафиксируем $0 \leq t_i \leq t_k \leq t_j \leq 1$ и $y, y_1, y_2 \in \mathbb{R}^d$. Для переносящих потоков $\theta^n_{t_i, t_j}(\cdot)$ и $\theta_{t_i, t_j}(\cdot)$ обозначим через $\theta^n_{t_j, t_i}(\cdot)$ и $\theta_{t_j, t_i}(\cdot)$ соответствующие прямые потоки. Тогда существует $C \geq 1$ такое, что:
		\item[\textbf{1)}]\begin{equation*}
				\left|\theta_{t_i, t_j}(y_1) - \theta_{t_i, t_j}(y_2)\right| + \left|\theta^n_{t_i, t_j}(y_1) - \theta^n_{t_i, t_j}(y_2)\right| \leq C \left|y_1 - y_2\right|;
		\end{equation*} 
		\item[\textbf{2)}]\begin{align*}
			\theta_{t_i, t_k}(\theta_{t_k, t_j}(y)) =  \theta_{t_i, t_j}(y), \\
			\theta^n_{t_i, t_k}(\theta^n_{t_k, t_j}(y)) =  \theta^n_{t_i, t_j}(y);
		\end{align*} 
		\item[\textbf{3)}]\begin{align*}
			C^{-1}\left|\theta_{t_j, t_i}(y_1) - y_2\right| \leq \left|\theta_{t_i, t_j}(y_2) - y_1\right| \leq C\left|\theta_{t_j, t_i}(y_1) - y_2\right|, \\
			C^{-1}\left|\theta^n_{t_j, t_i}(y_1) - y_2\right| \leq \left|\theta^n_{t_i, t_j}(y_2) - y_1\right| \leq C\left|\theta^n_{t_j, t_i}(y_1) - y_2\right|;
			\end{align*}
		\item[\textbf{4)}]\begin{equation*}
			C^{-1} \left(\left|y\right| \vee (t_j - t_i) \right)\leq \left|\theta_{t_i, t_j}(y)\right| + \left|\theta^n_{t_i, t_j}(y)\right| \leq C \left(\left|y\right| \vee (t_j - t_i) \right);
		\end{equation*} 
	\item[\textbf{5)}]\begin{equation*}
		\left|\theta_{t_i, t_j}(y) - \theta^n_{t_i, t_j}(y)\right| \leq C(t_j - t_i)\left(\Delta_b + \dfrac{1}{n^{\beta}}\right)(1 + |y|).
	\end{equation*} 
	\end{itemize}
\end{proposition}
\proof Доказательство свойства \textbf{1)} следует из цепочки неравенств 	\begin{align*}
	\left|\theta_{t, s}(x) - \theta_{t, s}(y)\right| = \left|x - y + \int_t^s \left(b(u, \theta_{u,s}(x)) - b(u, \theta_{u,s}(y))\right)du\right|  \leq \left|x - y\right| + C \cdot \int\limits_{t}^{s}\left|\theta_{u, s}(y) - \theta_{u, s}(x)\right|du
\end{align*} с последующим применением неравенства Гронуолла. Доказательство для дискретного потока аналогично. Свойство \textbf{2)} хорошо известно и немедленно влечет доказательство свойства \textbf{3)}. Далее, подставляя в утверждение из свойства \textbf{3)} $y_1 = 0$ и учитывая, что \begin{equation*}
\left|\theta_{t, s}(0)\right| = \left|\int_t^s \left(b(u, \theta_{u,s}(0)) \right)du\right| \leq C\left(s - t\right) + C  \int\limits_{t}^{s}\left|\theta_{u, s}(0)\right|du,
\end{equation*}
получаем требуемое. Докажем \textbf{5)}. Имеем: \begin{align*}
	\left|\theta^n_{t_i, t_j}(y) - \theta_{t_i, t_j}(y)\right| =  \left|\sum\limits_{k=i}^{j-1} \int_{t_k}^{t_{k+1}}  b_n(t_k, \theta^n_{t_k, t_j}(y)) - b(u, \theta_{u, t_j}(y)) du\right| \leq \\ \leq \left|\sum\limits_{k=i}^{j-1} \int_{t_k}^{t_{k+1}}  b_n(t_k, \theta^n_{t_k, t_j}(y)) - b_n(t_k, \theta_{t_k, t_j}(y)) du\right| +  \left|\sum\limits_{k=i}^{j-1} \int_{t_k}^{t_{k+1}}  b_n(t_k, \theta_{t_k, t_j}(y)) - b(t_k, \theta_{t_k, t_j}(y)) du\right| + \\ + \left|\sum\limits_{k=i}^{j-1} \int_{t_k}^{t_{k+1}}  b(t_k, \theta_{t_k, t_j}(y)) - b(t_k, \theta_{u, t_j}(y)) du\right| +  \left|\sum\limits_{k=i}^{j-1} \int_{t_k}^{t_{k+1}}  b(t_k, \theta_{u, t_j}(y)) - b(u, \theta_{u, t_j}(y)) du\right|. 
\end{align*}
Учитывая предположения \textbf{(A2)}, \textbf{(A3)} и \textbf{(A5)}, получаем:
\begin{align*}
	\left|\theta^n_{t_i, t_j}(y) - \theta_{t_i, t_j}(y)\right|   \leq \dfrac{C}{n}\sum\limits_{k=i}^{j-1}\left|   \theta^n_{t_k, t_j}(y) -\theta_{t_k, t_j}(y)\right| +  C(t_j - t_i)\Delta_{b} + \\ + \sum\limits_{k=i}^{j-1} \int_{t_k}^{t_{k+1}} \int_{t_k}^{u} \left|b(v, \theta_{v, t_j}(y))\right|dv du +  C\sum\limits_{k=i}^{j-1} \int_{t_k}^{t_{k+1}}  (u - t_k)^{\beta} du. 
\end{align*}
Применение неравенства Гронуолла доказывает утверждение \textbf{5)}.
\endproof

\section{Доказательство основного результата}

В данном разделе мы изучаем скорость (слабой) сходимости дискретного марковского процесса \eqref{МС_0} к диффузионному пределу \eqref{SDE_0}. Оба процесса допускают разложения в соответствующие ряды параметрикса \eqref{decomp_discr} и \eqref{decomp_cont}, что позволяет доказать основной результат работы \eqref{Main_Thm}, почленно сравнивая слагаемые указанных разложений, используя \begin{align}\label{decomp}
	&|p - p_n|(t_i, t_j, x, y) = |\sum_{r=0}^{\infty} \tilde{p} \otimes H^r(t_i, t_j, x, y) - \sum_{r=0}^{j-i} \tilde{p}^n \otimes_n H^{n,r}\left(t_i, t_j, x, y\right)| \leq \\ &\leq |\sum_{r=0}^{j-i} \tilde{p}^n \otimes_n (H^{n,r} - H^r) |\left(t_i, t_j, x, y\right) + |\sum_{r=0}^{j-i} (\tilde{p}^n - \tilde{p}) \otimes_n H^r |\left(t_i, t_j, x, y\right) + \\ &+ |\sum_{r=0}^{\infty} \tilde{p} \otimes_n H^r - \tilde{p} \otimes H^r |\left(t_i, t_j, x, y\right) =  i + ii + iii.
\end{align}
\begin{remark}
	При $r > j - i$ положим $\tilde{p} \otimes_n H^r = 0$. Доказательство Следствия 2.3 приводит к оценке 
	\begin{equation}
		|\sum_{r=j-i+1}^{\infty} \tilde{p} \otimes H^r|(t_i, t_j, x, y )\leq C (t_j - t_i)^{\gamma/2} \bar{p}(t_i, t_j, x, y).
	\end{equation} 
\end{remark}
Для оценки слагаемого $ii$ нам понадобится следующее утверждение:
\begin{lemma}
	Для любых $0 \leq t_i < t_j \leq 1$ и $x,y \in \mathbb{R}^d$, а также мультииндекса $\nu$ с $|\nu| \leq 4$ существует константа $C > 0$ такая, что справедливо 
	\begin{align}\label{diff_fr_dens}
		\left|D_x^{\nu}\tilde{p} -  D_x^{\nu}\tilde{p}_n\right|(t_i, t_j, x, y) &\leq \\ &\leq \frac{C(1 + |y|)}{(t_j - t_i)^{|\nu|/2}}\Delta^n\left(\mathcal{Q}_{S-d-2-|\nu|}\left(\frac{x - \theta_{t_i, t_j}(y)}{\sqrt{t_j - t_i}}\right) + \mathcal{Q}_{S-d-2-|\nu|}\left(\frac{x - \theta^n_{t_i, t_j}(y)}{\sqrt{t_j - t_i}}\right)\right).
	\end{align}
\end{lemma}
\proof Доказательство данного утверждения строится в два шага: на первом, аналогично доказательству дискретной части оценки \eqref{froz_dens_bound}, с помощью Теоремы 19.3 из \cite{BR} мы сравниваем переходное ядро $\tilde{p}_n$ с переходной плотностью $q$ вспомогательного гауссовского процесса, средние и ковариации которого совпадают со средними и ковариациями ведущего члена разложения в ряд Эджворта цепи \eqref{MC_1}. Имеем: 
\begin{align*}
	\left|D_x^{\nu}q -  D_x^{\nu}\tilde{p}_n\right|(t_i, t_j, x, y) \leq  \frac{C}{(t_j - t_i)^{|\nu|/2}}\dfrac{1}{\sqrt{n}} \mathcal{Q}_{S- d - 1 - |\nu|}\left(\frac{x - \theta^n_{t_i, t_j}(y)}{\sqrt{t_j - t_i}}\right).
\end{align*}
Второй шаг состоит в оценке разностей переходных плотностей $\tilde{p}$ и $q$ и их производных, то есть переходных плотностей гауссовских процессов с отличающимися средними и ковариациями. Для этого отождествим плотности $\tilde{p}$ и $q$ с $(d + 1) \times d$ матрицами $\Omega_{\tilde{p}}$ и $\Omega_q$, определенными при доказательстве Леммы 1 в \cite{KKM}. 
Далее, применяя разложение в ряд Тейлора первого порядка для функции 
\begin{align*}
	&f: \mathbb{R}^{(d + 1) \times d} \rightarrow \mathbb{R}, \\
	&A \mapsto f(A) =  \dfrac{1}{(2\pi)^{d/2}\operatorname{det}(A_{2 : d + 1})^{1/2}} &\exp \left(-\dfrac{1}{2} \langle (A_{2 : d + 1})^{-1}(y - A_1 - x), y - A_1 - x \rangle\right),
\end{align*} 
\begin{align*}
		\left|\tilde{p} - q\right|(t_i, t_j, x, y) = \left|f(\Omega_{\tilde{p}}) - f(\Omega_q)\right| =  \left| \sum\limits_{|\nu|=1}\left(\Omega_{\tilde{p}}-\Omega_q\right)^{\nu} \cdot\int\limits_{0}^{1}(1-\lambda)\mathcal{D}^{\nu}f\left\{\Omega_{\tilde{p}} + \lambda(\Omega_q-\Omega_{\tilde{p}})\right\}d\lambda \right|.
	\end{align*}
Для оценки правой части последнего равенства используем \eqref{froz_dens_bound} и оценки разностей средних и ковариаций. Для разности ковариационных матриц имеем:
\begin{align*}
	\left|\int_{t_i}^{t_j}a(u, \theta_{u, s}(y))du - \dfrac{1}{n}\sum_{k = i}^{j-1}a^n(t_k, \theta^n_{t_k, t_j}(y))\right| \leq \\ \leq \sum_{k = i}^{j-1} \int_{t_k}^{t_{k+1}}\left(\left|a(u, \theta_{u, s}(y)) - a^n(u, \theta_{u, s}(y))\right| + \left|a^n(u, \theta_{u, s}(y)) - a^n(t_k, \theta^n_{t_k, t_j}(y))\right|\right)du \leq \\ \leq C(t_j - t_i)\Delta_a + C(t_j - t_i) \cdot \left(\dfrac{1}{n^{\alpha}} + \Delta_{b}\right)(1 + |y|) \leq C(t_j - t_i)\Delta^n(1 + |y|).
\end{align*}

\begin{remark}
	Заметим, что ввиду очевидных неравенств 
	\begin{align}\label{change_sv}
		1 + |y| \leq 1 + |y - \theta^n_{t_j, t_i}(x)| + |\theta_{t_j, t_i}(x)| \leq  C\left(1 + |x -\theta^n_{t_i, t_j}(y)| + |\theta_{t_j, t_i}(x)|\right)
	\end{align}
	имеем
	\begin{align*}
		(1+|y|)\mathcal{Q}_{S - d - 2}\left(\frac{x -  \theta^n_{t_i, t_j}(y)}{\sqrt{t_j - t_i}}\right) \leq  (1+|x|)\mathcal{Q}_{S - d - 3}\left(\frac{x -  \theta^n_{t_i, t_j}(y)}{\sqrt{t_j - t_i}}\right).
	\end{align*}
\end{remark}
Таким образом, оценка \eqref{diff_fr_dens} также может быть записана в виде
\begin{equation*}
	\left|\tilde{p} - q\right|(t_i, t_j, x, y) \leq C(t_j - t_i)(1 + |x|)\Delta^n \mathcal{Q}_{S - d - 3}\left(\frac{x - \theta_{t_i, t_j}(y)}{\sqrt{t_j - t_i}}\right).
\end{equation*}
Это доказывает требуемое утверждение для случая $|\nu| = 0$. Остальные случаи рассматриваются аналогичным образом.
\endproof
Для дальнейшего анализа необходим аналог уравнения Колмогорова-Чепмэна для полиномиальных ядер $\mathcal{Q}_S(\dfrac{z}{\sqrt{t}})$.
\begin{proposition}
	Для всяких $0 \leq t_i < t_k < t_j \leq 1$ и $x,y \in \R^d$ существует $C \geq 1$ такое, что 
	\begin{equation}\label{conv_ker}
		\int_{\R^d} \mathcal{Q}_S\left(\dfrac{z - \theta_{t_i, t_k}(x)}{\sqrt{t_k - t_i}}\right) \mathcal{Q}_S\left(\dfrac{y - \theta_{t_k, t_j}(z)}{\sqrt{t_j - t_k}}\right)dz \leq C \mathcal{Q}_S\left(\dfrac{x -  \theta_{t_i, t_j}(y)}{\sqrt{t_j - t_i}}\right).
	\end{equation}
\end{proposition}
\proof См. \cite{KM0, KKM}.
\endproof
Теперь покажем, как связаны между собой полиномиальные ядра с одинаковым параметром при замене аргумента с $\frac{x - \theta^n_{t_i, t_j}(y)}{\sqrt{t_j - t_i}}$ на $\frac{x - \theta_{t_i, t_j}(y)}{\sqrt{t_j - t_i}}$. 
\begin{proposition}
	Для всяких $0 \leq t_i < t_j \leq 1$ и $x,y \in \R^d$ существует $C > 0$ такое, что 
	\begin{equation}\label{change}
\mathcal{Q}_{S}\left(\frac{x - \theta^n_{t_i, t_j}(y)}{\sqrt{t_j - t_i}}\right) \leq C(1 + |x|^S)\mathcal{Q}_{S}\left(\frac{x -  \theta_{t_i, t_j}(y)}{\sqrt{t_j - t_i}}\right).
	\end{equation}
	\end{proposition}
\proof Доказательство очевидным образом следует из неравенства 
\begin{equation*}
	\dfrac{1}{\left(1 + \dfrac{|x - \theta^n_{t_i, t_j}(y)|}{\sqrt{t_j - t_i}}\right)^S} \leq C \dfrac{\left(1 + \left|\theta^n_{t_i, t_j}(y) - \theta_{t_i, t_j}(y)\right|\right)^S}{\left(1 + \dfrac{|x - \theta_{t_i, t_j}(y)|}{\sqrt{t_j - t_i}}\right)^S}
\end{equation*}
и свойства \textbf{5)} переносящего потока.
\endproof
\begin{remark}
	Конечно, справедливо и неравенство
	\begin{equation*}
		\mathcal{Q}_{S}\left(\frac{x - \theta_{t_i, t_j}(y)}{\sqrt{t_j - t_i}}\right) \leq C(1 + |x|^S)\mathcal{Q}_{S}\left(\frac{x -  \theta^n_{t_i, t_j}(y)}{\sqrt{t_j - t_i}}\right).
	\end{equation*}
\end{remark}
Теперь мы готовы оценить слагаемое $ii$ в разложении \eqref{decomp}.
\begin{lemma}
	Для всяких $0 \leq t_i < t_j \leq 1$ и $x,y \in \R^d$ существует $C > 0$ такое, что
	\begin{equation}\label{2}
		ii \leq C(1 + |x|^{S - d - 2})\Delta^n \mathcal{Q}_{S - d - 3}\left(\frac{x -  \theta_{t_i, t_j}(y)}{\sqrt{t_j - t_i}}\right).
	\end{equation}
\end{lemma}
\proof Пусть $r = 1$. В этом случае из \eqref{ker_bound}, \eqref{diff_fr_dens}, \eqref{conv_ker} и \eqref{change} следует 
\begin{align*}
	&\left|(\tilde{p}_n - \tilde{p}) \otimes_n H\right|(t_i, t_j, x, y) \leq  \sum_{k = i}^{j - 1} \dfrac{1}{n}\left|\int_{\R^d} \left(\tilde{p}_n - \tilde{p}\right)(t_i, t_k, x, z) H(t_k, t_j, z, y)\right| dz \leq \\ &\leq C(1 + |x|^{S - d - 2})\Delta^n\mathcal{Q}_{S - d - 3}\left(\frac{x -  \theta_{t_i, t_j}(y)}{\sqrt{t_j - t_i}}\right) \int_{t_i}^{t_j} (t_j - u)^{\gamma/2 - 1} du = \\ &= C\dfrac{\Gamma(\gamma/2)}{\Gamma(1 + \gamma/2)}(t_j - t_i)^{\gamma/2}(1 + |x|^{S - d - 2})\Delta^n\mathcal{Q}_{S - d - 3}\left(\frac{x -  \theta_{t_i, t_j}(y)}{\sqrt{t_j - t_i}}\right).
\end{align*}
Предположим, что \begin{align*}
	&\left|(\tilde{p}_n - \tilde{p}) \otimes_n H^r\right|(t_i, t_j, x, y) \leq  C\dfrac{\Gamma^r(\gamma/2)}{\Gamma(1 + r\gamma/2)}(t_j - t_i)^{r\gamma/2}(1 + |x|^{S - d - 2})\Delta^n\mathcal{Q}_{S - d - 3}\left(\frac{x -  \theta_{t_i, t_j}(y)}{\sqrt{t_j - t_i}}\right).
\end{align*}
Тогда \begin{align*}
	&\left|(\tilde{p}_n - \tilde{p}) \otimes_n H^{r+1}\right|(t_i, t_j, x, y) \leq \\ &\leq C\dfrac{\Gamma^r(\gamma/2)}{\Gamma(1 + r\gamma/2)}(1 + |x|^{S - d - 2})\Delta^n\mathcal{Q}_{S - d - 3}\left(\frac{x -  \theta_{t_i, t_j}(y)}{\sqrt{t_j - t_i}}\right) \int_{t_i}^{t_j} (u - t_i)^{r\gamma/2} (t_j - u)^{\gamma/2 - 1} du = \\ &= C\dfrac{\Gamma^{r+1}(\gamma/2)}{\Gamma(1 + (r+1)\gamma/2)}(t_j - t_i)^{(r+1)\gamma/2}(1 + |x|^{S - d - 2})\Delta^n\mathcal{Q}_{S - d - 3}\left(\frac{x -  \theta_{t_i, t_j}(y)}{\sqrt{t_j - t_i}}\right).
\end{align*}
Это доказывает \eqref{2}.
\endproof

Перейдем теперь к оценке члена $iii$ разложения \eqref{decomp}. Обозначим $$p^{\otimes_n}(t_i, t_j, x, y) = \sum_{r=0}^{\infty} \left(\tilde{p} \otimes_n H^r\right)(t_i, t_j, x, y)$$. Для каждого слагаемого из $iii$ справедливо 
\begin{align*}
	\left(\tilde{p} \otimes H^r - \tilde{p} \otimes_n H^r\right)(t_i, t_j, x, y) &= \\ &= \left(\tilde{p} \otimes H^r - \tilde{p} \otimes H^{r-1} \otimes_n H\right)(t_i, t_j, x, y) + \\ &+ \left(\left[\tilde{p} \otimes H^{r-1} - \tilde{p} \otimes_n H^{r-1}\right] \otimes_n H\right)(t_i, t_j, x, y),
\end{align*}
что, в свою очередь, ведет к 
\begin{align*}
\left(p - p^{\otimes_n}\right)(t_i, t_j, x, y) = \sum_{r=0}^{\infty}	\left(\tilde{p} \otimes H^r - \tilde{p} \otimes_n H^r\right)(t_i, t_j, x, y)  = \\ = \left(p \otimes H - p \otimes_n H\right)(t_i, t_j, x, y) + \left[\left(p  - p^{\otimes_n}\right) \otimes_n H\right](t_i, t_j, x, y),
\end{align*}
а после итерирования данного соотношения получаем 
\begin{align}\label{iter}
	\left(p - p^{\otimes_n}\right)(t_i, t_j, x, y) = \\ = \left(p \otimes H - p \otimes_n H\right)(t_i, t_j, x, y) + \left[\left(p \otimes H - p \otimes_n H\right) \otimes_n \mathcal{H}\right](t_i, t_j, x, y),
\end{align}
где $\mathcal{H}(t_i, t_j, x, y) = \sum_{r=0}^{\infty} H^{r, \otimes_n}(t_i, t_j, x, y)$, а символ $H^{r, \otimes_n}$ обозначает $r$-кратное применение операции дискретной свертки к  ядру параметрикса $H$. Таким образом, ключевой становится оценка разности $\left(p \otimes H - p \otimes_n H\right)(t_i, t_j, x, y)$.

\begin{lemma}
	Для всяких $0 \leq t_i <  t_j \leq 1$ и $x,y \in \R^d$ существует $C > 0$ такое, что
	\begin{equation}\label{first_term_conv}
	\left|p \otimes H - p \otimes_n H\right|(t_i, t_j, x, y)	\leq C(1+ |x|^3)(t_j - t_i)^{\gamma/2}\bar{p}(t_i, t_j, x, y).
	\end{equation}
\end{lemma}
\proof \begin{align*}
	|\tilde{p} \otimes H - \tilde{p} \otimes_n H|(t_i, t_j, x, y)| = \\ = |\int_{t_i}^{t_j} \int_{\mathbb{R}^d} p(t_i, u, x, z)H(u, t_j, z, y) dz du - \frac{1}{n}\sum_{k=i}^{j-1} \int_{\mathbb{R}^d} p(t_i, t_k, x, z)H(t_k, t_j, z, y)dz| \leq \\ \sum_{k=i}^{j-1} \int_{t_k}^{t_{k+1}}\int_{\mathbb{R}^d} \left|\phi_u(z) - \phi_{t_k}(z)\right| dz du,
\end{align*}
где $\phi_u(z) = p(t_i, u, x, z)H(u, t_j, z, y)$. Разлагая подынтегральную разность в ряд Тейлора по переменной времени, получаем: \begin{equation*}
|\tilde{p} \otimes H - \tilde{p} \otimes_n H|(t_i, t_j, x, y)| \leq \sum_{k=i}^{j-1} \int_{t_k}^{t_{k+1}}\int_{\mathbb{R}^d} \int_0^1 (u - t_k) |\frac{d}{dv}\phi_v(z)|_{v = t_k + \delta(u - t_k)}d\delta dz du.
\end{equation*}
Далее, используя прямое и обратное равнения Колмогорова и сопряженность операторов в этих уравнениях, запишем  \begin{align*}
	\int_{\mathbb{R}^d} |\frac{d}{dv}\phi_v(z)|dz =  \int_{\mathbb{R}^d} p(t_i, v, x, z)(L_v^2 - 2L_v\tilde{L}_v + \tilde{L}_v^2)\tilde{p}(v, t_j, z, y) dz. 
\end{align*}

Обозначим $\psi_v(z) = L_v\tilde{p}(v, t_j, z, y)$. Аналогично $\tilde{\psi}_v(z) = \tilde{L}_v \tilde{p}(v, t_j, z, y)$. Таким образом, \begin{align*}
	&L_v\psi_v(z) - 2L_v\tilde{\psi}_v(z) + \tilde{L}_v\tilde{\psi}_v(z) =  \sum_{k=1}^d b_k(u,z) (\frac{d\psi_v}{dz_k} - \frac{d\tilde{\psi}_v}{dz_k}) + \\ &+ \sum_{k=1}^d (b_k(v,z) - b_k(v, \theta_{v, t_j}(y)) \frac{d\psi_v}{dz_k} + \dfrac{1}{2} \sum_{k, l=1} a_{kl}(v,z)(\frac{d^2\psi_v}{dz_kdz_l} - \frac{d\tilde{\psi}_v}{dz_kdz_l}) + \\ + &\dfrac{1}{2} \sum_{k, l=1} (a_{kl}(v,z) - a_{kl}(v, \theta_{v, t_j}(y))\frac{d\tilde{\psi}_v}{dz_kdz_l}.
\end{align*}
Теперь:
\begin{align*}
	&\left|\frac{d\psi_v}{dz_k}(z) - \frac{d\tilde{\psi}_v}{dz_k}(z)\right| \leq  \sum_{r=1}^d \left|\dfrac{d b_r(v,z)}{dz_k}\right|\left|\dfrac{d\tilde{p}}{dz_r}\right|  + \sum_{r=1}^d |b_r(v,z) - b_r(v, \theta_{v, t_j}(y))|\left|\dfrac{d^2\tilde{p}}{dz_k dz_r}\right| + \\ &+ \dfrac{1}{2} \sum_{r,l=1}^d \left|\dfrac{d a_{rl}(v,z)}{dz_k}\right|\left|\dfrac{d^2\tilde{p}}{dz_r dz_l}\right| + \dfrac{1}{2} \sum_{r,l=1}^d |a_{rl}(v,z) - a_{rl}(v, \theta_{v, t_j}(y))|\left|\dfrac{d^3\tilde{p}}{dz_k dz_r dz_l}\right| \leq \\ &\leq C \bar{p}(v, t_j, z, y) \left(\dfrac{1 + |z|}{\sqrt{t_j - v}} + \dfrac{|z - \theta_{v, t_j}(y)|}{t_j - v} + \dfrac{1}{t_j - v} + \dfrac{|z - \theta_{v, t_j}(y)|^{\gamma}}{(t_j - v)^{3/2}}\right) \leq \\ &\leq \dfrac{C}{(t_j - v)^{(3-\gamma)/2}}(1 + |z|)\bar{p}(v, t_j, z, y).
\end{align*}
Аналогично получаем, что 
\begin{equation*}
	\left|\frac{d^2\psi_v}{dz_k dz_l}(z) - \frac{d^2\tilde{\psi}_v}{dz_k dz_l}(z)\right| \leq \dfrac{C}{(t_j - v)^{2-\gamma/2}}(1 + |z|)\bar{p}(v, t_j, z, y).
\end{equation*}
Далее, применяя полученные оценки, имеем

\begin{align*}
	|\tilde{p} \otimes H - \tilde{p} \otimes_n H|(t_i, t_j, x, y)| \leq \\ \leq C(1+|x|^3)\bar{p}(t_i, t_j, x, y)\sum_{k=i}^{j-1} \int_{t_k}^{t_{k+1}} \int_0^1 (u - t_k)(t_j - t_k - \delta(u - t_k))^{\gamma/2 - 2}d\delta du = \\ =  C(1+|x|^3)\bar{p}(t_i, t_j, x, y)\sum_{k=i}^{j-1} \int_{t_k}^{t_{k+1}}  \left((t_j-u)^{\gamma/2-1} - (t_j - t_k)^{\gamma/2 - 1}\right) du \leq \\ \leq C(1+|x|^3)\bar{p}(t_i, t_j, x, y)\int_{t_i}^{t_{j}} (t_j-u)^{\gamma/2-1} du = C(t_j - t_i)^{\gamma/2}(1+|x|^3)\bar{p}(t_i, t_j, x, y).
\end{align*}
\endproof

Учитывая оценку \begin{equation*}
	\left|\mathcal{H}\right|(t_i, t_j, x, y) \leq \dfrac{C}{(t_j - t_i)^{1-\gamma/2}}\bar{p}(t_i, t_j, x, y),
\end{equation*}
\eqref{first_term_conv} и применяя уравнение Колмогорова-Чепмэна для переходной плотности $\bar{p}$ вспомогательного диффузионного процесса, получаем из \eqref{iter}
\begin{equation}\label{3}
	iii \leq C(t_j - t_i)^{\gamma/2}(1 + |x|^3)\bar{p}(t_i, t_j, x, y).
\end{equation}

Наибольшую   представляет оценка слагаемого $i$ в разложении \eqref{decomp}. Проведем ее в 2 этапа. На первом этапе мы почленно сравним сумму $\sum_{r=0}^{j-i} \tilde{p}_n \otimes_n H^r$ с ${\sum_{r=0}^{j-i} \tilde{p}_n \otimes_n \H^r}$,
где \begin{equation*}
	\H\left(t_i, t_j, x, y\right):=\left(L_{t_i}-\tilde{L}_{t_i}\right) \tilde{p}_n\left(t_{i}, t_j, x, y\right).
\end{equation*}

\begin{lemma}
	Для всяких $0 \leq t_i <  t_j \leq 1$ и $x,y \in \R^d$ существует $C > 0$ такое, что
	\begin{align}\label{fir_part}
		\sum_{r=0}^{j-i}\left|\tilde{p}_n \otimes_n (H^r - \H^r)\right|(t_i, t_j, x, y) &\leq \\ &\leq C(1 + |x|^{S-d-3} + |y|^{S-d-3})\Delta^n\mathcal{Q}_{S - d - 5}\left(\frac{x -  \theta^n_{t_i, t_j}(y)}{\sqrt{t_j - t_i}}\right).
	\end{align}
\end{lemma}
\proof  Для случая $r = 1$  имеем 

\begin{align*}
	\left|H - \H\right|(t_i, t_j, x, y) = \left| \left(L_{t_i}-\tilde{L}_{t_i}\right)\left(\tilde{p} - \tilde{p}_n\right)\right|(t_i, t_j, x, y) \leq \\ \leq \sum_{k=1}^d \left|b_k(t_i, x) - b_k(t_i, \theta_{t_i, t_j}(y))\right|\left|\dfrac{d(\tilde{p} - \tilde{p}_n)}{dx_k}\right| + \dfrac{1}{2} \sum_{k, l=1}^d \left|a_{kl}(t_i, x) - a_{kl}(t_i, \theta_{t_i, t_j}(y))\right|\left|\dfrac{d^2(\tilde{p} - \tilde{p}_n)}{dx_kdx_l}\right|. 
\end{align*}

Далее, из линейности роста переносящих потоков, \eqref{diff_fr_dens}, \eqref{froz_dens_bound}, Замечания \eqref{change_sv} и очевидной цепочки неравенств 
\begin{align*}
	1 + |y| \leq C \left(1 + |z - \theta_{t_i, t_j}(y)| + |x - \theta^n_{t_i, t_j}(z)| + |x|\right)
\end{align*}
следует, что
\begin{align*}
	\left|\tilde{p}_n \otimes_n (H - \H)\right|(t_i, t_j, x, y) \leq  \dfrac{C\Gamma(\gamma/2)}{\Gamma(1 + \gamma/2)}(t_j - t_i)^{\gamma/2}(1 + |x|^{S - d - 3})\Delta^n\mathcal{Q}_{S - d - 5}\left(\frac{x -  \theta^n_{t_i, t_j}(y)}{\sqrt{t_j - t_i}}\right).
\end{align*}

Далее, 
\begin{align*}
	\left|\tilde{p}_n \otimes_n \H\right|(t_i, t_j, x, y) \leq  \dfrac{C\Gamma(\gamma/2)}{\Gamma(1 + \gamma/2)}{(t_j - t_i)^{\gamma/2}}\mathcal{Q}_{S - d - 5}\left(\frac{x -  \theta^n_{t_i, t_j}(y)}{\sqrt{t_j - t_i}}\right),
\end{align*}
и, аналогично оценке \eqref{conv_convergence}, по индукции
\begin{align*}
	\left|\tilde{p}_n \otimes_n \H^r\right|(t_i, t_j, x, y) \leq  \dfrac{C\Gamma^r(\gamma/2)}{\Gamma(1 + r\gamma/2)}(t_j - t_i)^{r\gamma/2}\mathcal{Q}_{S - d - 5}\left(\frac{x -  \theta^n_{t_i, t_j}(y)}{\sqrt{t_j - t_i}}\right).
\end{align*}
Теперь, итерируя равенство
\begin{align*}
	\tilde{p}_n \otimes_n H^{r + 1} - \tilde{p}_n \otimes_n \H^{r + 1} =  \tilde{p}_n \otimes_n \left( H^{r} - \H^{r}\right) \otimes_n H  + \left(\tilde{p}_n \otimes_n \H^{r} \right) \otimes_n \left(\H - H\right),
\end{align*}
получаем оценку 
\begin{align*}
	\left|\tilde{p}_n \otimes_n H^{r + 1} - \tilde{p}_n \otimes_n \H^{r + 1}\right| &\leq \\ &\leq \left|\tilde{p}_n \otimes_n (H - \H) \otimes_n H^{\otimes_n, r}\right|  + \sum_{k=0}^{r-1}\left| \left(\tilde{p}_n \otimes_n \H^{r-k} \right) \otimes_n \left(\H - H\right) \otimes_n H^{\otimes_n, k}\right|.
\end{align*}
Заметим, что для функции $f(u) = (u-t_i)^{r\gamma/2}(t_j-u)^{\gamma/2-1}$ ввиду монотонного роста $f$ на $[t_i, t_j]$ при $r \geq 0$ справедливо 
\begin{equation*}
	\sum_{k=0}^{j-i-1}\dfrac{1}{n}f(t_{i+k}) \leq \int_{t_i}^{t_j} f(u)du.
\end{equation*} 
Подставляя полученные выше оценки и применяя \eqref{conv_ker}, имеем 
\begin{align*}
	\left|\tilde{p}_n \otimes_n (H - \H) \otimes_n H^{\otimes_n, r}\right| \leq \\ \leq\dfrac{C\Gamma^{r+1}(\gamma/2)}{\Gamma(1 + (r+1)\gamma/2)}(t_j - t_i)^{(r+1)\gamma/2}(1 + |x|^{S - d - 3})\Delta^n\mathcal{Q}_{S - d - 5}\left(\frac{x -  \theta^n_{t_i, t_j}(y)}{\sqrt{t_j - t_i}}\right); \\
	\left| \left(\tilde{p}_n \otimes_n \H^{r-k} \right) \otimes_n \left(\H - H\right) \otimes_n H^{\otimes_n, k}\right| \leq \\ \leq \dfrac{C\Gamma^{r+1}(\gamma/2)}{\Gamma(1 + (r+1)\gamma/2)}(t_j - t_i)^{(r+1)\gamma/2}(1 + |y|^{S - d - 3})\Delta^n\mathcal{Q}_{S - d - 5}\left(\frac{x -  \theta^n_{t_i, t_j}(y)}{\sqrt{t_j - t_i}}\right).
\end{align*}
Доказательство утверждения \eqref{fir_part} завершается после применения формулы Стирлинга для асимптотического поведения Гамма-функции.
\endproof

Перейдем теперь к сравнению ядра $\H$ с дискретным ядром параметрикса.

\begin{lemma}
	Для всяких $0 \leq t_i <  t_j \leq 1$ и $x,y \in \R^d$ существует $C > 0$ такое, что
	\begin{align}\label{sec_part}
		\sum_{r=0}^{j-i}\left|\tilde{p}_n \otimes_n (H_n^r - \H^r)\right|(t_i, t_j, x, y) &\leq \\ &\leq C\ln(e(j-i))(1 + |x|^{S-d-2})\Delta^n\mathcal{Q}_{S - d - 6}\left(\frac{x -  \theta^n_{t_i, t_j}(y)}{\sqrt{t_j - t_i}}\right).
	\end{align}
\end{lemma}
\proof
Применяя разложения Тейлора 3-го порядка для дискретных инфинитезимальных генераторов, примененных к плотности $\tilde{p}_n$

\begin{align*}
	H_n(t_i, t_j, x, y) = \left({L}_{t_i}^n - \tilde{L}_{t_i}^n\right)\tilde{p}_n(t_i, t_j, x, y) = \\ = n\int_{\R^d} q_{t_i, x}^n(z)\left(\tilde{p}_n(t_i, t_j, x + \dfrac{1}{n}b_n(t_i, x) + \dfrac{1}{\sqrt{n}}z, y) - \tilde{p}_n(t_i, t_j, x, y)\right)dz - \\ - n\int_{\R^d} q_{t_i, \theta^n_{t_i, t_j}(y)}^n(z)\left(\tilde{p}_n(t_i, t_j, x + \dfrac{1}{n}b_n(t_i, \theta^n_{t_i, t_j}(y)) + \dfrac{1}{\sqrt{n}}z, y) - \tilde{p}_n(t_i, t_j, x, y)\right)dz = \\ = n\int_{\R^d} q_{t_i, x}^n(z)\left(\sum_{k=1}^d \left(\dfrac{1}{n}(b_n)_k(t_i, x) + \dfrac{1}{\sqrt{n}}z_k\right)\dfrac{d\tilde{p}_n}{dx_k}(t_i, t_j, x, y) + \right. \\  \left. + \dfrac{1}{2}\sum_{k,l=1}^d \left(\dfrac{1}{n}(b_n)_k(t_i, x) + \dfrac{1}{\sqrt{n}}z_k\right)\left(\dfrac{1}{n}(b_n)_l(t_i, x) + \dfrac{1}{\sqrt{n}}z_l\right)\dfrac{d^2\tilde{p}_n}{dx_kdx_l}(t_i, t_j, x, y) + \right. \\ \left.  + 3\sum_{|\nu|=3}\dfrac{1}{\nu !}\left(\dfrac{1}{n}b_n(t_i, x) + \dfrac{1}{\sqrt{n}}z\right)^{\nu}\int_0^1 (1-\delta)^2 D_x^{\nu}\tilde{p}_n(t_i, t_j, x + \dfrac{\delta}{n}b_n(t_i, x) + \dfrac{\delta}{\sqrt{n}}z, y)d\delta\right) dz - \\ - n\int_{\R^d} q_{t_i, \theta^n_{t_i, t_j}(y)}^n(z)\left(\sum_{k=1}^d \left(\dfrac{1}{n}(b_n)_k(t_i, \theta^n_{t_i, t_j}(y)) + \dfrac{1}{\sqrt{n}}z_k\right)\dfrac{d\tilde{p}_n}{dx_k}(t_i, t_j, x, y) + \right. \\  \left. + \dfrac{1}{2}\sum_{k,l=1}^d \left(\dfrac{1}{n}(b_n)_k(t_i, \theta^n_{t_i, t_j}(y)) + \dfrac{1}{\sqrt{n}}z_k\right)\left(\dfrac{1}{n}(b_n)_l(t_i, \theta^n_{t_i, t_j}(y)) + \dfrac{1}{\sqrt{n}}z_l\right)\dfrac{d^2\tilde{p}_n}{dx_kdx_l}(t_i, t_j, x, y) + \right. \\ \left.  + 3\sum_{|\nu|=3}\dfrac{1}{\nu !}\left(\dfrac{1}{n}b_n(t_i, \theta^n_{t_i, t_j}(y)) + \dfrac{1}{\sqrt{n}}z\right)^{\nu}\int_0^1 (1-\delta)^2 D_x^{\nu}\tilde{p}_n(t_i, t_j, x + \dfrac{\delta}{n}b_n(t_i, \theta^n_{t_i, t_j}(y)) + \dfrac{\delta}{\sqrt{n}}z, y)d\delta\right) dz.
\end{align*}

С учётом данного разложения и свойств семейства плотностей $q^n_{t_i, x}(\cdot)$

\begin{align}
	&\left(H_n - \H\right)(t_i, t_j, x, y) = \\ &= \sum_{k=1}^d \left((b_n)_k(t_i, x) - (b_n)_k(t_i, \theta^n_{t_i, t_j}(y)) - b_k(t_i, x) + b_k(t_i, \theta_{t_i, t_j}(y)) \right)\dfrac{d\tilde{p}_n}{dx_k}(t_i, t_j, x, y) + \\ &+ \dfrac{1}{2} \sum_{k,l=1}^d \left(a^n_{kl}(t_i, x) - a^n_{kl}(t_i, \theta^n_{t_i, t_j}(y)) - a_{kl}(t_i, x) + a_{kl}(t_i, \theta_{t_i, t_j}(y)) \right)\dfrac{d^2\tilde{p}_n}{dx_kdx_l}(t_i, t_j, x, y) + \\ &+ \dfrac{1}{2n} \sum_{k,l=1}^d \left((b_n)_k(t_i, x)(b_n)_l(t_i, x) - (b_n)_k(t_i, \theta^n_{t_i, t_j}(y))(b_n)_l(t_i, \theta^n_{t_i, t_j}(y))\right)\dfrac{d^2\tilde{p}_n}{dx_kdx_l}(t_i, t_j, x, y) + \\ &+ \dfrac{3}{\sqrt{n}} \sum_{|\nu|=3}\dfrac{1}{\nu !}\int_{\R^d}\left(q^n_{t_i,x}(z)\left(\left(\dfrac{1}{\sqrt{n}}b_n(t_i, x) + z\right)^{\nu}\int_0^1 (1-\delta)^2 D_x^{\nu}\tilde{p}_n(t_i, t_j, x + \dfrac{\delta}{n}b_n(t_i, x) + \dfrac{\delta}{\sqrt{n}}z, y)d\delta\right) - \right. \\ &\left. - q^n_{t_i,\theta^n_{t_i, t_j}(y)}(z)\left(\left(\dfrac{1}{\sqrt{n}}b_n(t_i, \theta^n_{t_i, t_j}(y)) + z\right)^{\nu}\int_0^1 (1-\delta)^2 D_x^{\nu}\tilde{p}_n(t_i, t_j, x + \dfrac{\delta}{n}b_n(t_i, \theta^n_{t_i, t_j}(y)) + \dfrac{\delta}{\sqrt{n}}z, y)d\delta\right)\right) dz = \\ &= S_1 + S_2 + S_3 + S_4.
\end{align}

Учитывая предположение \textbf{(A2)}, \eqref{froz_dens_bound}, свойство \textbf{5)} переносящих потоков и Замечание \eqref{change_sv}, для $S_1$ и $S_2$ справедливо 

\begin{align}\label{S_1-S_2}
	|S_1| \leq 	\sum_{k=1}^d \left(\left|(b_n)_k(t_i, x) - b_k(t_i, x)\right| + \left|(b_n)_k(t_i, \theta^n_{t_i, t_j}(y)) - b_k(t_i, \theta^n_{t_i, t_j}(y))\right| + \right. \\ \left. + \left|b_k(t_i, \theta^n_{t_i, t_j}(y)) - b_k(t_i, \theta_{t_i, t_j}(y))\right|\right)|\dfrac{d\tilde{p}_n}{dx_k}(t_i, t_j, x, y)| \leq \\ \leq \dfrac{C}{\sqrt{t_j - t_i}}(1+|y|)\Delta^n \mathcal{Q}_{S - d - 3}\left(\frac{x -  \theta^n_{t_i, t_j}(y)}{\sqrt{t_j - t_i}}\right); \\
	|S_2| \leq \dfrac{1}{2}\sum_{k,l=1}^d \left(\left|a^n_{kl}(t_i, x) - a_{kl}(t_i, x)\right| + \left|a^n_{kl}(t_i, \theta^n_{t_i, t_j}(y)) - a_{kl}(t_i, \theta^n_{t_i, t_j}(y))\right| + \right.\\ \left. + \left|a_{kl}(t_i, \theta^n_{t_i, t_j}(y)) - a_{kl}(t_i, \theta_{t_i, t_j}(y))\right|\right)|\dfrac{d^2\tilde{p}_n}{dx_kdx_l}(t_i, t_j, x, y)| \leq \\ \leq \dfrac{C}{t_j - t_i}(1+|y|^{\gamma})\Delta^n \mathcal{Q}_{S - d - 4}\left(\frac{x -  \theta^n_{t_i, t_j}(y)}{\sqrt{t_j - t_i}}\right).
\end{align} 

Также, предположение \textbf{(A3)} приводит к 
\begin{align*}
	|S_3| \leq \dfrac{C}{t_j - t_i}(1+|x|^2+|y|^2)\Delta^n\mathcal{Q}_{S - d - 4}\left(\frac{x -  \theta^n_{t_i, t_j}(y)}{\sqrt{t_j - t_i}}\right) \leq \dfrac{C}{t_j - t_i}(1+|y|^2)\Delta^n\mathcal{Q}_{S - d - 6}\left(\frac{x -  \theta^n_{t_i, t_j}(y)}{\sqrt{t_j - t_i}}\right).
\end{align*}

Теперь зафиксируем некоторый мультииндекс $\eta$ с $|\eta|=3$ и введем обозначения:
\begin{align*}
	\rho_1(z) = q^n_{t_i,x}(z), \rho_2(z) = q^n_{t_i,\theta^n_{t_i, t_j}(y)}(z); \\
	\omega_1(z) = \int_0^1 (1-\delta)^2 D_x^{\eta}\tilde{p}_n(t_i, t_j, x + \dfrac{\delta}{n}b_n(t_i, x) + \dfrac{\delta}{\sqrt{n}}z, y)d\delta, \\ \omega_2(z) = \int_0^1 (1-\delta)^2 D_x^{\eta}\tilde{p}_n(t_i, t_j, x + \dfrac{\delta}{n}b_n(t_i, \theta^n_{t_i, t_j}(y)) + \dfrac{\delta}{\sqrt{n}}z, y)d\delta; \\
	\varepsilon_1(z) = \left(\dfrac{1}{\sqrt{n}}b_n(t_i, x) + z\right)^{\eta}, \varepsilon_2(z) = \left(\dfrac{1}{\sqrt{n}}b_n(t_i, \theta^n_{t_i, t_j}(y)) + z\right)^{\eta}.
\end{align*}

Тогда оценивание $S_4$ сведется к рассмотрению 
\begin{align*}
	&\int_{\mathbb{R}^d}|\omega_1(z)\varepsilon_1(z)\rho_1(z) - \omega_2(z)\varepsilon_2(z)\rho_2(z)|dz \leq \\  &\leq  \int_{\R^d}\left(|(\omega_1(z) - \omega_2(z))\varepsilon_1(z)\rho_1(z)| +  |\omega_2(z)(\varepsilon_1(z) - \varepsilon_2(z))\rho_1(z)| +  |\omega_2(z)\varepsilon_2(z)(\rho_1(z) - \rho_2(z))|\right)dz = \\ &= 1) + 2) + 3).
\end{align*}

Из предположения \textbf{(A4)}, оценки \eqref{froz_dens_bound}

\begin{align*}
	&3) \leq \dfrac{C|x - \theta^n_{t_i, t_j}(y)|}{\sqrt{n}}\int_{\R^d}(1 + |y| + |z|)^3Q_S(z) \times \\ &\times\int_0^1 (1-\delta)^2\left| D_x^{\eta}\tilde{p}_n(t_i, t_j, x + \dfrac{\delta}{n}b_n(t_i, \theta^n_{t_i, t_j}(y)) + \dfrac{\delta}{\sqrt{n}}z, y)\right|d\delta dz \leq  \dfrac{C|x - \theta^n_{t_i, t_j}(y)|}{(t_j - t_i)^{3/2}} \mathcal{Q}_{S-d-4}\left(\frac{x -  \theta^n_{t_i, t_j}(y)}{\sqrt{t_j - t_i}}\right) \times\\ &\times \int_{\R^d} \dfrac{(1 + |y| + |z|)^3}{(1 +|z|^S)}\int_0^1 (1-\delta)^2 \left(1 + \dfrac{\delta|b_n(t_i, \theta^n_{t_i, t_j}(y))|}{n} + \dfrac{\delta|z|}{\sqrt{n}}\right)^{S-d-4}d\delta dz \leq \\ &\leq \dfrac{C}{t_j - t_i}(1+ |y|^{S-d-1})\Delta^n\mathcal{Q}_{S-d-5}\left(\frac{x -  \theta^n_{t_i, t_j}(y)}{\sqrt{t_j - t_i}}\right).
\end{align*}

Далее, из разложения 
\begin{align*}
	&a_1 a_2 a_3 - b_1 b_2 b_3 = \\ = &(a_1 - b_1)a_2 a_3 + b_1 (a_2 - b_2) a_3 + b_1 b_2 (a_3 - b_3)
\end{align*}

аналогичным образом следует 
\begin{align*}
	2)&\leq \\ &\leq \dfrac{C|x - \theta^n_{t_i, t_j}(y)|}{\sqrt{n}}\int_{\R^d}Q_S(z)(1 + |y| + |z|)^2  \int_0^1 (1-\delta)^2\left| D_x^{\eta}\tilde{p}_n(t_i, t_j, x + \dfrac{\delta}{n}b_n(t_i, \theta^n_{t_i, t_j}(y)) + \dfrac{\delta}{\sqrt{n}}z, y)\right|d\delta dz \leq   \\ &\leq \dfrac{C}{t_j - t_i}(1+ |y|^{S-d-2})\Delta^n\mathcal{Q}_{S-d-6}\left(\frac{x -  \theta^n_{t_i, t_j}(y)}{\sqrt{t_j - t_i}}\right).
\end{align*}
Для оценки $1)$ применим разложение в ряд Тейлора первого порядка, а также предположение \textbf{(A2)} и оценку \eqref{froz_dens_bound}:

\begin{align*}
	&1) \leq C\int_{\R^d} q^n_{t_i, x}(z)(1 + |y|+|z|)^3 \times \\ &\times \int_0^1 (1-\delta)^2\left| D_x^{\eta}\tilde{p}_n(t_i, t_j, x + \dfrac{\delta}{n}b_n(t_i, x) + \dfrac{\delta}{\sqrt{n}}z, y) - D_x^{\eta}\tilde{p}_n(t_i, t_j, x + \dfrac{\delta}{n}b_n(t_i, \theta^n_{t_i, t_j}(y)) + \dfrac{\delta}{\sqrt{n}}z, y)\right|d\delta dz \leq \\ &\leq C\int_{\R^d} \dfrac{(1 + |y|+|z|)^3}{1 + |z|^S} \int_{0}^{1} \dfrac{\delta(1-\delta)^2}{n}\left|b_n(t_i, x) - b_n(t_i, \theta^n_{t_i, t_j}(y))|\right| \times \\ &\times \int_0^1 \sum_{k=1}^d \left|D_x^{\eta+e_k}\tilde{p}_n(t_i, t_j, x + \dfrac{\delta}{\sqrt{n}}z + \dfrac{(1 -\tau)\delta}{n}b_n(t_i, x) + \dfrac{\tau\delta}{n}b_n(t_i, \theta^n_{t_i, t_j}(y)), y)\right|d\tau d\delta dz \\ &\leq \dfrac{C}{(t_j-t_i)^2}\int_{\R^d} \dfrac{(1 + |y|+|z|)^3}{1 + |z|^S} \int_{0}^{1} \dfrac{\delta(1-\delta)^2}{n}\left|x - \theta^n_{t_i, t_j}(y)\right| \times \\ &\times \int_0^1 \sum_{k=1}^d \mathcal{Q}_{S-d-5}\left(\frac{x + \dfrac{\delta}{\sqrt{n}}z + \dfrac{\delta}{n}b_n(t_i, x) + \dfrac{\tau\delta}{n}\left(b_n(t_i, \theta^n_{t_i, t_j}(y)) - b_n(t_i, x)\right) - \theta^n_{t_i, t_j}(y)}{\sqrt{t_j - t_i}}\right) d\tau d\delta dz \leq \\ &\leq \dfrac{C}{n(t_j-t_i)^2}\int_{\R^d} \dfrac{(1 +|y|+|z|)^3}{1 + |z|^S} \int_{0}^{1} \delta(1-\delta)^2\left|x -  \theta^n_{t_i, t_j}(y)\right| \times \\ &\times (1 + |y| + |z|)^{S-d-5}\mathcal{Q}_{S-d-5}\left(\frac{x  - \theta^n_{t_i, t_j}(y)}{\sqrt{t_j - t_i}}\right) d\delta dz \leq \dfrac{C}{t_j - t_i}(1+ |y|^{S-d-2})\Delta^n\mathcal{Q}_{S-d-6}\left(\frac{x -  \theta^n_{t_i, t_j}(y)}{\sqrt{t_j - t_i}}\right).
\end{align*}

Последнее неравенство справедливо, поскольку $\sqrt{n} (t_j - t_i)^{1/2} \geq 1$.  

Таким образом, \begin{align*}
	|S_4| \leq \dfrac{C}{t_j - t_i}(1+ |y|^{S-d-2})\Delta^n\mathcal{Q}_{S-d-6}\left(\frac{x -  \theta^n_{t_i, t_j}(y)}{\sqrt{t_j - t_i}}\right),
\end{align*}

и, как следствие, 

\begin{align*}
	\left|H_n - \H\right|(t_i, t_j, x, y) \leq \dfrac{C}{t_j - t_i}(1+ |y|^{S-d-2})\Delta^n\mathcal{Q}_{S-d-6}\left(\frac{x -  \theta^n_{t_i, t_j}(y)}{\sqrt{t_j - t_i}}\right).
\end{align*} 
Учитывая \eqref{change_sv} и повторяя проведенные вычисления, получаем, для $k = 1\cdots 4$ мы можем записать полученные оценки в виде 
\begin{align*}
|S_k| \leq \dfrac{C}{t_j - t_i}(1+ |x|^{S-d-2})\Delta^n\mathcal{Q}_{S-d-6}\left(\frac{x -  \theta^n_{t_i, t_j}(y)}{\sqrt{t_j - t_i}}\right).
\end{align*}
Этот факт будет использован при оценке сверток высокого порядка.

Для дискретной свертки плотности "замороженного" процесса $\tilde{p}_n$ с разностью ядер $H_n - \H$ первого порядка справедливо 
\begin{align*}
	&\left|\tilde{p}_n \otimes_n (H_n - \H)\right|(t_i, t_j, x, y) \leq  \sum_{k=i}^{j-1} \dfrac{1}{n}\int_{\R^d} |\tilde{p}_n|(t_i, t_k, x, z)|H_n - \H|(t_k, t_j, z, y)dz \leq \\ &\leq C(1+|y|^{S-d-2})\Delta^n\mathcal{Q}_{S-d-6}\left(\frac{x -  \theta^n_{t_i, t_j}(y)}{\sqrt{t_j - t_i}}\right)\sum_{k=i}^{j-1} \dfrac{1}{n}\dfrac{1}{t_j - t_k} \leq \\ &\leq C(1+|y|^{S-d-2})\Delta^n\mathcal{Q}_{S-d-6}\left(\frac{x -  \theta^n_{t_i, t_j}(y)}{\sqrt{t_j - t_i}}\right) \sum_{k=1}^{j-i} \dfrac{1}{k} \leq \\ &\leq C \ln(e(j-i))(1+|y|^{S-d-2})\Delta^n\mathcal{Q}_{S-d-6}\left(\frac{x -  \theta^n_{t_i, t_j}(y)}{\sqrt{t_j - t_i}}\right).
\end{align*}

Далее, 
\begin{align*}
	\left|\left(\tilde{p}_n \otimes_n \H^{r} \right) \otimes_n \left(H_n - \H\right) \right|(t_i, t_j, x, y) \leq \\ \leq \dfrac{C^r\Gamma^r(\gamma/2)}{\Gamma(1 + r\gamma/2)}(1 +  |y|^{S-d-2})\mathcal{Q}_{S - d - 6}\left(\frac{x -  \theta^n_{t_i, t_j}(y)}{\sqrt{t_j - t_i}}\right)\Delta^n \sum_{k = i}^{j-1}\dfrac{1}{n}\dfrac{(t_k - t_i)^{r\gamma/2}}{t_j - t_k} \leq \\ \leq \dfrac{C^r\Gamma^r(\gamma/2)}{\Gamma(1 + r\gamma/2)}\ln(e(j-i))(1+|y
	|^{S-d-2})\Delta^n\mathcal{Q}_{S-d-6}\left(\frac{x -  \theta^n_{t_i, t_j}(y)}{\sqrt{t_j - t_i}}\right).
\end{align*}
Теперь, применяя итеративный спуск, как при доказательстве утверждения \eqref{fir_part}, получаем оценку \eqref{sec_part}.
\endproof
	
	Таким образом,
\begin{align}\label{1}
	i \leq C\ln(e(j-i))(1 + |x|^{S-d-2} + |y|^{S-d-2})\Delta^n\mathcal{Q}_{S-d-6}\left(\frac{x -  \theta^n_{t_i, t_j}(y)}{\sqrt{t_j - t_i}}\right),
\end{align}
что после сложения с оценками \eqref{change} \eqref{2} и \eqref{3} доказывает утверждение \eqref{Main_Thm}.

}


\begin{thebibliography}{9}
\bibitem{B25}Биттер И.И. Локальная предельная теорема для возмущенных выборочных траекторий индуцированных порядковых статистик // Управление большими системами. - 2025. - Вып. 113. - С.6-20.
\bibitem{Ar1}D. G. Aronson, The fundamental solution of a linear parabolic equation containing a small
parameter. {\em Ill. Journ. Math.} {\bf 3}(1959), 580--619.
\bibitem{Ar2}D. G. Aronson, Bounds for the fundamental solution of a parabolic equation. {\em Bull. Amer.
	Math. Soc.} {\bf 73}(1967), 890--896.
\bibitem{bass}R. F. Bass, Diffusions and Elliptic Operators. {\em Springer}, (1997).
\bibitem{BK21} Bitter, I. and Konakov, V. L1 and $L_{\infty}$ stability of transition densities of perturbed
diffusions. Random Oper. Stoch. Eq. 2021, 29 (4), 287 - 308.
\bibitem{BR} Bhattacharya, R. and Rao, R. Normal approximations and asymptotic expansions. Wiley
and sons, 1976.
\bibitem{David} David, H.A. (1973), Concomitants of order statistics, Bull. Internat.
Statist. Inst. 45, 295-300.
\bibitem{David2} David, H.A. and Galambos, J. (1974). The asymptotic theory of concomitants of order statistics, J. Appl. Probab. 11, 762-770.
\bibitem{Dav} Davydov, Yu. and Egorov, V. Functional limit theorems for induced order statistics. Mathematical Methods of Statistics, 9(3):297-313, (January, 2000).
\bibitem{deck}T. Deck, S. Kruse, Parabolic differential equations with Holder continuous and unbounded
coefficients. {\em Acta Applicandae Mathematicae 74} {\bf 1}(2002), 71--91.
\bibitem{DM}F.Delarue, S.Menozzi, Density estimates for a random noise propagating through a
chain of differential equations. {\em J. Funct. Anal. 259} {\bf 6}(2010), 1577--1630.
\bibitem{fried}A. Friedman, Partial Differential Equations of Parabolic Type. Prentice-Hall, (1964).
\bibitem{kalash}A.M. Il'in, A.S. Kalashnikov and O.A. Oleinik, Second-order linear equations of parabolic
type. {\em Uspehi Mat. Nauk} {\bf 17}(1962), 3--146.
\bibitem{KM0} Konakov, V. and Mammen, E. Local Limit Theorems and Strong Approximations for
Robbins-Monro Procedures. arXiv:2304.10673, (2023).
\bibitem{KM} Konakov, V. and Mammen, E. Local limit theorems for transition densities of Markov chains converging to diffusions.  Probability Theory and Related Fields, Volume 117, pages 551-587, (2000).
\bibitem{KKM}V.Konakov, A. Kozhina and S. Menozzi, Stability of densities for perturbed diffusions and
Markov chains. {\em ESAIM: Probability and Statistics} {\bf 21}(2017), 88--112.
\bibitem{KMM}V. Konakov, S. Menozzi and S. Molchanov, Explicit parametrix and local limit theorems
for some degenerate diffusion processes. {\em Annales de l'Institut Henri Poincare} {\bf 46}(2010),
908--923.
\bibitem{MPZ}S. Menozzi, A. Pesce, X. Zhang, Density and gradient estimates for non-degenerate
Brownian SDEs with unbounded measurable drift. {\em J. Diff. Eq.} {\bf 272}(2021), 330--369.
\bibitem{Sk1}Skorohod, A.V. Asymptotic methods for Stochastic differential equations. (in russian). Kiev
Naukova dumka (1987).
\bibitem{Sk2} Skorohod, A.V. Studies in the theory of random processes. Addison-Wesley, Reading, Massachussetts.[English translation of Skorohod, A. V. (1961). Issledovaniya po teorii sluchaynykh protsessov. Kiev University Press] (1965).
\bibitem{SV79} Stroock, D.W. and Varadhan, S.R.  Multidimensional diffusion processes.
Springer, Berlin, Heidelberg, New York (1979).
\end{thebibliography}
\end{document}